# BEYOND PAIR CORRELATION

Hugh L. Montgomery[1]
K. Soundararajan[2]

*Dedicated to the memory of Erdős Pál.*

Goldston & Montgomery [3] showed that under the assumption of the Riemann Hypothesis (RH), the Pair Correlation Conjecture of Montgomery [5] is equivalent to the assertion that

$$\int_1^X (\psi(x+h)-\psi(x)-h)^2\,dx \sim hX\log\frac{X}{h} \tag{1}$$

for $X^\epsilon \le h \le X^{1-\epsilon}$. In contrast, the Cramér model, which holds that the primes are distributed as if the integer $n$ is prime with probability $1/\log n$, each one independent of another, would predict that this expression is $\sim hX\log X$. If the Cramér model does not apply, one is left to speculate about the distribution of $\psi(x+h)-\psi(x)$. Recently the authors [6] used a quantitative form of the Prime $k$-tuple Hypothesis to give a heuristic determination of the moments of $\psi(x+h)-\psi(x)-h$, which supports the notion that $\psi(x+h)-\psi(x)$ is approximately normally distributed with mean $\sim h$ and variance $\sim h\log X/h$, as $x$ varies, $1\le x\le X$, with $h$ in the range $X^\epsilon\le h\le X^{1-\epsilon}$. Odlyzko [7] and Forrester & Odlyzko [2] analyzed the distribution of the zeros of the zeta function, and found that the data is in close agreement with the Pair Correlation Conjecture. Hence one might expect that numerical studies of primes in short intervals would lend support to the conjectural relation (1). With this in mind we have calculated the distribution of $\psi(x+h)-\psi(x)-h$ for $0\le x\le X=10^{10}$ when $h=10^5$. In Table 1 below we give the numerical values of the moments

$$\mu_k(X,h)=\frac{1}{X}\int_0^X(\psi(x+h)-\psi(x)-h)^k\,dx,$$

as well as of the normalized moments $\widetilde{\mu}_k=\mu_k/\mu_2^{k/2}$. Since the normal distribution has normalized moments $\widetilde{\mu}_{2k+1}=0$, $\widetilde{\mu}_{2k}=(2k-1)\cdot(2k-3)\cdots3\cdot1$, we see that the normalized moments are reasonably close to their anticipated values. The sixth moment is a little large, which suggests that large deviations may be rather more common than would otherwise be the case. In this regard we note that the largest value of $\psi(x+h)-\psi(x)-h$ encountered is

[1]Supported in part by a grant from the Number Theory Foundation.
[2]Supported in part by a grant from the American Institute of Mathematics.

5046.08 at $x = 9559758537$, which is 5.30 times the standard deviation. In $10^5$ independent samples, which is essentially what we presume to have here, the likelihood of such a large deviation occurring is $1 - \Phi(5.3)^{10^5} = 0.00577$. Here $\Phi(x) = \frac{1}{\sqrt{2\pi}}\int_{-\infty}^{x} e^{-t^2/2}\,dt$ is the cumulative distribution function of the normal variable with mean 0 and standard deviation 1. Similarly, the smallest value found is $-4920.06$ at $x = 5116809527$. This is $-5.17$ times the standard deviation; such a large negative value would occur, in $10^5$ independent samples of a normal variable, with probability $1 - \Phi(5.17)^{10^5} = 0.01163$. These large deviations are somewhat larger than might be expected, but not *so* much larger, since the maximum is larger than 4138 with probability 1/2. Finally, it was found that

$$\operatorname{meas}\{x \in [0, 10^{10}] : |\psi(x + 10^5) - \psi(x) - 10^5| > 3000\} = 3080882.$$

Since the size of this set is less than one fifth the size one would expect with a comparable normal variable, the large deviations at this threshhold are less common than would be predicted.

| $k$ | $\mu_k$ | $\widetilde{\mu}_k$ |
|---|---|---|
| 0 | $1.0000$ | 1.0000 |
| 1 | $9.0984 \times 10^{-2}$ | 0.0001 |
| 2 | $9.0663 \times 10^{5}$ | 1.0000 |
| 3 | $-1.1926 \times 10^{6}$ | −0.0014 |
| 4 | $2.4995 \times 10^{12}$ | 3.0408 |
| 5 | $-2.4951 \times 10^{13}$ | −0.0319 |
| 6 | $1.1573 \times 10^{19}$ | 15.5288 |

Table 1. Moments of $\psi(x + h) - \psi(x) - h$ for $0 \le x \le X = 10^{10}$ with $h = 10^5$.

In addition to the numerical data described above, the results of sieving were also recorded in the form of the cumulative distribution function, and plotted against that of a normal variable with the same variance, in Figure 1. The fit to normal is impressive. Note that both functions are being graphed on the same coordinate axes.

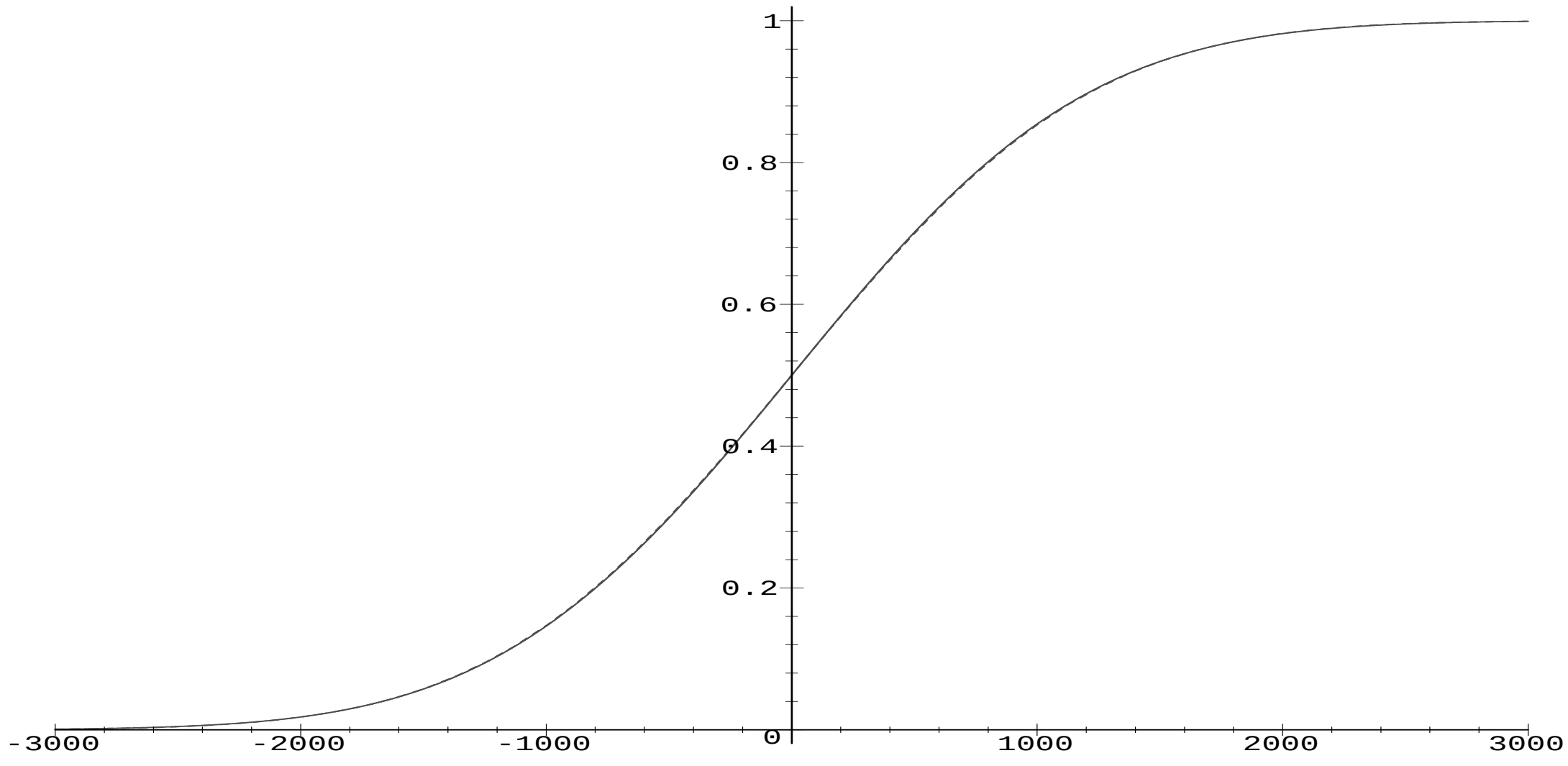

Figure 1. Distribution of $\psi(x + h) - \psi(x) - h$ (solid) versus normal (dashed).

One of the objects of the numerical study was to test whether the variance of $\psi(x+h)-\psi(x)-h$ is near the value $h\log X = 23.02\times 10^5$ that would be predicted by the Cramér model, or whether it is nearer the to the smaller variance $h\log X/h = 11.51\times 10^5$ predicted by (1). The big surprise in the data is that the variance $9.07\times 10^5$ recorded in Table 1 is significantly smaller than even the smaller of these values. To address this discrepancy we reconsider the heuristics used to develop (1). Upon expanding, we see that the left hand side of (1) is approximately

$$\sum_{m\le X}\sum_{n\le X}\Lambda(m)\Lambda(n)\max(0,h-|m-n|)-h^2X.$$

This in turn is approximately

$$h\sum_{n\le X}\Lambda(n)^2+2\sum_{k=1}^{h}(h-k)\sum_{n\le X}\Lambda(n)\Lambda(n+k)-h^2X.$$

By using the Prime Number Theorem with a sharp remainder (we may assume RH), we see that the first term above is approximately $hX\log X-hX$. As for the second term, we let $E(X,k)$ be defined by the relation

$$\sum_{n\le X}\Lambda(n)\Lambda(n+k)=\mathfrak{S}(k)X+E(X,k)$$

where $\mathfrak{S}(k)$ is the singular series defined by Hardy & Littlewood [4] for the Twin Prime Conjecture,

$$\mathfrak{S}(k)=\prod_{p|k}\Big(1+\frac{1}{p-1}\Big)\prod_{p\nmid k}\Big(1-\frac{1}{(p-1)^2}\Big).$$

If $k$ is odd then $\mathfrak{S}(k)=0$, but if $k$ is even then

$$\mathfrak{S}(k)=c\prod_{\substack{p|k\\ p>2}}\frac{p-1}{p-2}$$

where

$$c=2\prod_{p>2}\Big(1-\frac{1}{(p-1)^2)}\Big).$$

It is well-known that $\mathfrak{S}(k)$ is 1 on average, and the estimate with Cesàro weights,

$$\sum_{k=1}^{h}(h-k)\mathfrak{S}(k)=\frac{1}{2}h^2-\frac{1}{2}h\log h+O(h)$$

was used by Montgomery (1971, unpublished) to guess at the Pair Correlation Conjecture. We now refine this estimate.

**Theorem.** *Let $\mathfrak{S}(k)$ be defined as above. Then*

$$\sum_{k=1}^{h}(h-k)\mathfrak{S}(k) = \frac{1}{2}h^2 - \frac{1}{2}h\log h + Ah + O(h^{1/2+\epsilon})$$

*where $A = (1 - C_0 - \log 2\pi)/2$. (Here $C_0$ is Euler's constant.)*

When we insert this in the earlier calculation, we come to the conclusion that we should expect that

$$(2) \qquad \int_0^X \big(\psi(x+h) - \psi(x) - h\big)^2\, dx = hX\log\frac{X}{h} + BhX + \text{smaller terms}$$

where $B = -C_0 - \log 2\pi = -2.41509\ldots$. For $X = 10^{10}$ and $h = 10^5$, this more accurate main term predicts a second moment of $9.098 \times 10^5$, which is much closer to the computed value, $9.066 \times 10^5$.

The main barrier to majorizing the 'smaller terms' in (2) lies in estimating the contribution

$$2\sum_{k=1}^{h}(h-k)E(X,k)$$

of the error terms in the Twin Prime Conjecture. Numerical studies (cf. Brent [1]) suggest that $E(X,k) \ll X^{1/2+\epsilon}$, and one may presume that this holds uniformly for $1 \le k \le X$. Thus the above quantity should be $\ll h^2X^{1/2+\epsilon}$, but we actually expect that there is some cancellation in the sum itself, so that the above is $\ll h^{3/2+\epsilon}X^{1/2+\epsilon}$. Indeed, when all the possible sources of error are taken into account, one concludes that the relation (2) may hold with an error term that is $\ll h^{1/2}X^{1/2+\epsilon} + h^{3/2+\epsilon}X^{1/2}$.

*Proof of the Theorem.* Let $s(k) = \prod_{p|k, p>2}\frac{p-1}{p-2}$. Then

$$\sum_{k=1}^{h}(h-k)\mathfrak{S}(k) = c\sum_{k=1}^{h/2}(h-2k)s(2k) = 2c\sum_{k=1}^{h/2}(h/2-k)s(k).$$

We show that

$$\sum_{k=1}^{K}(K-k)s(k) = \frac{K^2}{c} - \frac{K\log K}{2c} + \frac{K}{2c}(1 - C_0 - \log 4\pi),$$

which suffices. Let

$$S(s) = \sum_{k=1}^{\infty} s(k)k^{-s} = \big(1-2^{-s}\big)^{-1}\prod_{p>2}\Big(1 + \frac{p-1}{(p-2)(p^s-1)}\Big)$$

for $\Re s > 1$. Then

$$S(s) = \zeta(s)\prod_{p>2}\Big(1 + \frac{1}{(p-2)p^s}\Big) = \zeta(s)T(s),$$

say, for $\Re s > 0$. Similarly, we note that

$$T(s) = \zeta(s+1)\big(1-2^{-s-1}\big)\prod_{p>2}\Big(1+\frac{2}{(p-2)p^{s+1}}-\frac{1}{(p-2)p^{2s+1}}\Big) = \zeta(s+1)\big(1-2^{-s-1}\big)U(s),$$

say, for $\Re s > -1/2$. Clearly,

$$\sum_{k=1}^{K}(K-k)s(k) = \frac{1}{2\pi i}\int_{a-i\infty}^{a+i\infty} S(s)\frac{K^{s+1}}{s(s+1)}\,ds$$

when $a$ is a real number, $a > 1$. We move the integral to the abscissa $b$, where $-1/2 < b < 0$, and consider the residues arising from the simple pole in the integrand at $s = 1$ and the double pole at $s = 0$. Since $\zeta(s) \sim 1/(s-1)$ when $s$ is near 1, and since $T(1) = 2/c$, it follows that the residue at $s = 1$ is $K^2/c$. As for the residue at $s = 0$, we recall from Titchmarsh [8, pp. 16–20] that

$$\zeta(s+1) = \frac{1}{s} + C_0 + O(|s|), \qquad \zeta(0) = -1/2, \qquad \zeta'(0) = -\frac{1}{2}\log 2\pi.$$

Also, $U(0) = 2/c$ and $U'(0) = 0$. Hence, with a little calculation, we see that the residue at $s = 0$ is

$$-\frac{K\log K}{2c} + \frac{K}{2c}(1 - C_0 - \log 4\pi).$$

As for the remaining integral, we note by the functional equation and Stirling's formula that $|\zeta(b+it)| \asymp V^{1/2-b}$ when $V \le t \le 2V$. Also, by the Cauchy–Schwarz inequality,

$$\int_V^{2V} |\zeta(b+1+it)|\,dt \le V^{1/2}\Big(\int_V^{2V} |\zeta(b+1+it)|^2\,dt\Big)^{1/2} \ll_b V,$$

in view of known mean-square estimates of the zeta function (cf. Theorem 7.2 of Titchmarsh [8]). Since $U(b+it) \ll_b 1$ for $b > -1/2$, it follows that the integral in question is absolutely convergent with a value $\ll_b K^{b+1}$. Since we may take $b$ as close to $-1/2$ as we please, this gives the stated result.

When approached as above, it seems fortuitous that $T(1) = U(0) = 2/c$ and that $U'(0) = 0$. But miracles do not happen by accident, so it seems that there is something going on here that remains to be understood.

## References

1. R. P. Brent, *Irregularities in the distribution of primes and twin primes*, Math. Comp. **29** (1975), 43–56; Correction **30** (1976), 198.
2. P. J. Forrester & A. M. Odlyzko, *Gaussian unitary ensemble eigenvalues and Riemann $\zeta$ function zeros: a nonlinear equation for a new statistic*, Phys. Rev. E (3) **54** (1996), R4493–R4495.
3. D. A. Goldston & H. L. Montgomery, *On pair correlations of zeros and primes in short intervals*, Analytic Number Theory and Diophantine Problems (Stillwater, OK, July 1984) (A. C. Adolphson, J. B. Conrey, A. Ghosh, R. I. Yager, eds.), Prog. Math. 70, Birkäuser, Boston, 1987, pp. 183–203.

4. G. H. Hardy & J. E. Littlewood, *Some problems of "Partitio Numerorum" (III): On the expression of a number as a sum of primes*, Acta Math. **44** (1922), 1–70.
5. H. L. Montgomery, *The pair correlation of zeros of the zeta function*, Analytic number theory (St. Louis Univ., 1972), Proc. Sympos. Pure Math. **24**, Amer. Math. Soc., Providence, 1973, pp. 181–193.
6. H. L. Montgomery & K. Soundararajan, *Primes in short intervals*, to appear.
7. A. M. Odlyzko, *On the distribution of spacings between zeros of the zeta function*, Math. Comp. **48** (1987), 273–308.
8. E. C. Titchmarsh, *The Theory of the Riemann Zeta Function*, Second Edition, Oxford University Press, Oxford, 1986.

Dept. of Math., Univ. of Michigan, Ann Arbor, MI 48109–1009, USA
email: hlm@math.lsa.umich.edu

Institute for Advanced Study, Princeton, NJ 08540, USA
email: ksound@math.ias.edu